\documentstyle[numreferences]{kluwer}
\input amssym.def
\input amssym

\begin{opening}
\title{Covariant $q$-differential Calculus \\
and its Deformations at $q^N =1$}
\author{ Richard {\surname Kerner}}
\author{ Bernd \surname{Niemeyer}\thanks{
Stagiaire DEA, \'Ecole Normale Sup\'erieure de Lyon, 
46 All\'ee d'Italie, Lyon, France.
Present address: Department of Physics, University of Hannover,
Germany}}
\institute{Laboratoire de Gravitation et Cosmologie Relativistes\\
 Universit\'e Pierre-et-Marie-Curie - CNRS URA D0769\\
 Tour 22, 4$^{eme}$ \'etage, Bo\^{\i}te 142,\\
 4, Place Jussieu, 75005 Paris, FRANCE}
\date{February 21,1998}
\end{opening} 
\begin{document}
\begin{abstract}  {We construct the generalized version of covariant 
${\Bbb Z}_3$-graded differential calculus introduced by one of us (R.K.) 
in  \cite{Ker1}, \cite{Ker2}, and then extended to the case of arbitrary 
${\Bbb Z}_N$ grading in (\cite{Kapr1}, \cite{MDVRK1}, \cite{MDVRK2}). 
Here our main purpose is to establish the recurrence formulae for the $N$-th 
power of covariant $q$-differential $D_q = d_q + A$ and to analyze more 
closely the particular case of $q$ being an $N^{\rm {th}}$ primitive root 
of unity. The generalized notions of connection and curvature are introduced 
and several examples of realization are displayed for $N=3$ and $4$. Finally 
we briefly discuss the idea of infinitesimal deformations of the parameter 
$q$ in the complex plane.}
\end{abstract}

\section {Introduction}

The idea of investigating the generalizations of exterior differential
calculus by postulating $d^N = 0$ instead of the usual $d^2 = 0$, leading 
to generalized Grassmann algebras with $n$-th order constitutive relations is 
not really a new one, and has been mentioned here and there quite a long time 
ago, but only recently it has been given the attention it really deserves.
In a recent series of articles (\cite{Kapr1}, \cite{Ker1}, \cite{Ker2}, 
\cite{Ker3}, \cite{MDVRK1}, \cite{MDVRK2}, \cite{MDV1}), the case $N = 3$ has
been investigated in more detail (cf. also \cite{Ker4}, \cite{AKLR1}), then
a general theory of $q$-differentials with $q = N$ and $d^N = 0$ has been
worked out, introducing also the $q$-analogs of the homological algebras.
\newline
\indent
Besides the universal construction in terms of tensor products of algebras
and linear spaces, some simple concrete realizations have been found, namely
a finite (matrix) version in which the operator $d$ is identified with a
${\Bbb Z}_N$ commutator with a grade $1$ element of the ${\Bbb Z}_N$-graded 
algebra, and a differential calculus on the infinite space of functions 
generated by a finite number of variables $x^k$ has been also defined 
(\cite{Ker1}, \cite{Ker2}). Finally, the notion of the {\it covariant 
differential} has been discussed, and generalized for the ${\Bbb Z}_3$-graded 
case. In particular, it enabled us to introduce the generalized notion of the 
{\it curvature form}, which in the ${\Bbb Z}_3$-graded case was equal to the 
$3$-form $D^2 A$, the second covariant differential of the connection 
$1$-form $A$.
\newline
\indent
In the present article we shall investigate the general $q$-deformed 
differential calculus, with special attention being paid to the case when $q$
is equal to an $N$-th primitive root of unity, $q^N = 1$. We shall show how
in this case the $N$-th covariant differential acting on an appropriate 
module ${\cal{H}}$ reduces itself to an automorphism of ${\cal{H}}$, i.e. 
$D^N \Phi = \Omega \, \Phi$, with $\Omega = D^{N-1} A$. We shall also derive
a simple recurrence formula for $D^N \, \Phi$ for an arbitrary value of the 
parameter $q$.
\newline
\indent
Next, we discuss briefly two realizations of this $q$-differential calculus
with $q$ an $N$-th primitive root of unity: tha $N \times N$ complex matrix 
representation, and a generalized Grassmannian spanned by all the differential
forms of the type $ d x^k , \, d^2 x^k , \, \dots d^{N-1} x^k $. The explicit
expressions are found for $\Omega \, = \, D^{N-1} A $ in the matrix case, and
for a few low values of $N$ in the generalized Grassmannian case. In the
case of matrix realization, the general form of a matrix representing the
``pure gauge'' configuration is given.
\newline
\indent
Then we consider the {\it infinitesimal deformations} of the complex parameter 
$q$ itself. It seems worthwhile to know what happens when the consecutive 
infinitesimal deformations form a polygon in ${\Bbb C}^1$, so that after $N$ 
steps we come back to the initial value of $q$. In the case when $q$ tends to 
the $N$-th primitive root of unity, certain combinations of products $N$ 
deformed differentials yield $0$-degree operators with interesting properties.
\newline
\indent
Finally, we look at the infinitesimal deformations of the covariant 
$q$-differential $D_q$, in which not only the first term $d_q$ is transformed
into $d_{q+\epsilon}$, but parallelly also the connection one-form $A$ is
replaced by $A + \epsilon \, \Lambda$. Here again, we compute the $N$-th 
order products of the deformed covariant differentials.

 \section{Universal $q$-differential and its properties}

Let ${\cal{A}}$ be an associative algebra with unit element, generated by two 
elements denoted by $U$ and $\eta$, respectively. We attribute the
{\it grade $1$} to the element $\eta$, and {\it grade $0$} to the element $U$:
 $ {\rm {deg}}(U)=  0 $, ${\rm {deg}}(\eta) = 1$.\\
$q$ being a complex number different from $0$, we shall impose the following 
commutation relation between the generators $U$ and $\eta$:
\begin{equation}
U \, \eta = q \, \eta \, U
\end{equation}
A general element belonging to ${\cal{A}}$ can now be represented by a finite
sum of various powers of the generators $U$ and $\eta$ in ordered products:
\begin{equation}
B \in {\cal{A}}, \, \ \ \, B = {\displaystyle{\sum_{m,n}}} \, b_{mn} \, U^m \, 
\eta^n = {\displaystyle{\sum_n}}\, \beta_n(U)\,\eta^n
\end{equation}
\indent
An element of {\cal{A}} has a well-defined degree  $n$ if it is a monomial
of $n$-th order in the generator $\eta$. The algebra ${\cal{A}}$ acquires
a natural ${\Bbb Z}$-grading and can be represented as an infinite sum of 
subspaces with well-defined grades:
$$ {\cal{A}} \, = \, {\cal{A}}_0 \oplus {\cal{A}}_1 \oplus {\cal{A}}_2 \oplus
{\cal{A}}_3 \oplus \dots \oplus {\cal{A}}_k \oplus {\cal{A}}_{k+1} \oplus \dots
$$
\indent
If $q$ is a primitive $N$-th root of unity, $q^N = 1$, and if the following 
supplementary conditions are imposed on the generators $U$ and $\eta$ :
\begin{equation}
U^N \, = \, {\bf 1}, \, \ \ \, \eta^N \, = \, {\bf 1} ,
\end{equation}
then our algebra becomes finite and can be represented as algebra of complex
$N \times N$ matrices. Its canonical representation will be introduced later
in the Section 4. In that case the matrix $U$ is called the {\it grading 
matrix}, the matrix $\eta$ has an inverse which is equal to $\eta^T$, and 
Eq. (1) can be written as:
$$U^{-1} \, \eta \, U \, = \, q \, \eta , $$
consequently,
\begin{equation}
U^{-1} \, \eta^m \, U = q^m \, \eta^m \label{commU}
\end{equation}
\indent
The {\it $q$-differential} of an element $B$ of grade $b$ in ${\cal{A}}$ can 
be defined now as follows (for any value of the parameter $q$:)
\begin{equation}
d_q \, B = \eta \, B - q^b \, B \, \eta 
\end{equation}
The grades add up under the associative multiplication in our algebra; that
is why the $q$-differential satisfies the $q$-deformed Leibniz rule. For
a product of two matrices $A \, B$, with deg$(A) = a$ and deg$(B)=b$ one has
\begin{equation}
d_q \, (A \, B) \, = \, (d_q \, A) \, B \, + \, q^a \, A \, (d_q B)
\end{equation}
The proof is straightforward by applying the definition:

$$ d_q (A \, B) \, = \, \eta \, (AB) - q^{a+b} \, (AB) \, \eta \, =$$
$$\, \eta \, AB \, - \, q^a \, A \eta B + q^a \, A \eta B -
q^a q^b \, AB \, \eta =$$
$$  \, [ \eta \, A - q^a \, A\, \eta ] \, B + q^a \, A \, [ \eta \, B - q^b 
\, B \, \eta] \, = (d_q \, A) \, B \, + \, q^a \, A \, (d_q B)  $$
\newline
\indent
Now we evaluate the consecutive powers of the operator $d_q$:
$$d_q^2 \, B \, = 
d_q \, ( \eta B - q^b \, B \eta) \, = \, \eta ( \eta B - q^b \, B \eta ) - 
q^{b+1} \, ( \eta B - q^b \, B \eta ) \, \eta \, =  $$ 
$$ \, \eta^2 \, B - q^b \, (1 + q) \, \eta B \eta \, + q^{2 b + 1} \,
B \, \eta^2 ; $$
$$ d_q^3 \, B \, = \, \eta^3 \, B + q^b (1+q+q^2) \eta^2 B \eta + q^{2 b + 1} 
(1 + q + q^2) \, \eta B \eta^2 - q^{3 b + 3} \, B \eta^3 $$
\indent
It is enough to check the action of consecutive powers of the operator $d_q$
on the two generators $U$ and $\eta$ in order to be able to extend them on
an arbitrary element $B$  ($b = {\rm {deg}}(B)$) of the entire algebra 
${\cal{A}}$. We find out easily that for the element of degre one, $\eta$, 
we have
$$d_q \, \eta = \eta^2 - q \, \eta^2 ,$$
$$d^2_q \, \eta = (1-q) \, \eta^3 - (1-q) \, q^2 \, \eta^3 =
(1-q)(1-q^2) \, \eta^3,..., $$
\begin{equation}
d^N_q \, \eta = (1-q)(1-q^2) \dots (1-q^N) \, \eta^{N+1}
\end{equation}
If the condition $q^N = 1$ is imposed, the last expression vanishes, implying 
$d^N \eta = 0.$ When the operator $d_q$ acts of the $0$-degree generator $U$, 
we find
$$d_q \, U = \eta \, U - U \, \eta, $$
$$d^2_q \, U = \eta^2 \, U - \eta U \eta - q \, \eta U \eta + q \, U \eta^2$$
Using the $q$-deformed commutation relation \ref{commU} we can write
$$- \eta U \eta = - q \, \eta^2 \, U \, \ \ \, \ \ {\rm and} \, \ \ \, \ \
- q \, \eta U \eta = - U \, \eta^2 $$
so that one gets 
$$d^2_q \, U = (1-q) \, (\eta^2 \, U - U \, \eta^2). $$
Similarly, one has for $d^N \, U :$
\begin{equation}
d^N_q \, U = (1-q)(1-q^2)\dots(1-q^{N-1}) \, (\eta^N \, U - U \, \eta^N)
\end{equation}
which vanishes because we have assumed $\eta^N = {\bf 1 }$, so that it does
commute with any element of the algebra.
\newline
\indent

It is easy to see that independently of the grade of $B$, one has 
$d^N \, B = 0$ when $q$ is the $N$-th primitive root of unity and when
$\eta^N = {\bf 1}$. This is obvious in the two particular cases shown above,
with $q = - 1$ and $q = j = e^{2 \pi i / 3}$. The general formula can be 
easily proved by recurrence using the properties of the generalized binomial 
symbols 
$${N \brack k}_q = \frac{[N]_q ! }{ [k]_q ! [(N-k)]_q !}, $$
\begin{equation}
{\rm with} 
\, \ \ \, \ \ [k]_q = (1 + q + q^2 + ... + q^k) \, \ \  {\rm  and \,} \, \ \ 
[k]_q ! = [1]_q \, [2]_q \, [3]_q ... [k]_q ,    \label{factorials1}
\end{equation}
\indent
Classical recurrence relations remain valid, their proof by induction obvious:
\begin{equation} 
 [n+1]_q \, [n]_q! = [n+1]_q !  \ \ \, {\rm \, and \, } 
\ \ \, \ \  {\displaystyle{n\brack k}_q}  \, = \,
\frac{[n]_q!}{[n-k]_q! \, [k]_q!}   \label{factorials2}
\end{equation}
\begin{equation} 
 \, \ \ {\rm and} \, \ \ \, q^{k+1}\,{n\brack{k+1}}_q \,+ \, 
{n\brack k}_q \,= \,{{n+1}\brack{k+1}}_q   \label{factorials3}
\end{equation}
\indent
The general formula for $d_q^N \, B$ can be also established by induction and
reads as follows:
\begin{equation}
d^N_q \, B = {\displaystyle{\sum}_{k = 0}^N} \, (-1)^k \, q^{k \, b +
\frac{k(k-1)}{2}} \, {\displaystyle{N \brack k}_q} \, \eta^{N - k} \, B \, 
\eta^k
\end{equation}
\indent
Now we can easily check that $d^N_q \, B = 0$. Indeed, the $q$-binomial 
coefficients $\displaystyle{N \brack{k}}$ vanish for $0 < k < N$ because of
the common factor $[N]_q = (1 + q + q^2 +...+ q^{N-1} = 0$, the only remaining 
term being equal to
$$d^N_q \, B = \eta^N \, B + (-1)^N \, q^{N b + \frac{N(N-1)}{2}} \, B \, 
\eta^N$$
Now, because $\eta^N={\bf 1}$ is the unit element of the algebra ${\cal{A}}$,
commuting with any $B \in {\cal{A}}$, and because of $q^{N b} = 1$, we can
write
$$ d_q^{N} \, B = \biggl( 1 + (-1)^N \, q^{\frac{N(N-1)}{2}} \biggr) \, B $$
This expression vanishes identically whether $N$ is even or odd number. If $N$
is odd, then $(N-1)/2 = c$ is an integer, so that $q^{c N} = 1$ , while
$(-1)^N = -1$. If $N$ is even, then $q^{\frac{N}{2}} = - 1$, $(-1)^N = 1$,
but $(q^{\frac{N}{2}})^{(N-1)}) = (-1)^{(N-1)} = - 1$, so in both cases we
have a factor $1 + (-1)$ in front of $B$, which makes this expression vanish,
thus completing the proof that $d^N_q \, B = 0$

\section{ Covariant $q$-differential and its successive powers}

Defining the covariant $q$-differential acting on an element of the module
$\Phi \in {\cal{H}}$ as above, i.e.
$$D \, \Phi = d \Phi \, + \, A\, \Phi ,$$
$A$ denoting the connection $1$-form, which is an element of degree $1$
belonging to the algebra ${\cal{A}}$, and using the $q$-Leibniz rule, we get
the following explicit expressions for the consecutive powers of $D$ acting
on $\Phi$, i.e. the formulae for $D^N \, \Phi$:
\begin{equation}
D \Phi = d \Phi + A\,\Phi;
\end{equation}
\begin{equation}
D^2 \Phi = d^2 \Phi + (1+q)\, A\,d \Phi + (D A)\,\Phi ;
\end{equation}
In the limit of $q = -1$ we have $d^2 \Phi = 0$, so that only the last term
of the above expression survives, yielding the well-known formula for 
curvature in the ${\Bbb Z}_2$-graded case, $D^2 \Phi = (DA)\,\Phi$. Now,
\begin{equation}
D^3 \Phi = d^3 \Phi + (1+q+q^2)\, A \,d^2 \Phi + (1+q+q^2) \, (D A) \, 
d^2 \Phi + (D^2 A) \, \Phi;
\end{equation}
Here again, when $q = j = e^{2 \pi i/3}$ is the primitive $3$-rd root of
unity, we have $d^3 \Phi = 0$ by definition, and the coefficient $(1+q+q^2)$
vanishes, leaving only $D^3 \Phi = (D^2 A)\, \Phi$.  Similarly, 
$D^4 \Phi = (D^3 A)\, \Phi $, and so on. 
\newline
\indent
The general formula can be quite easily established with the help of notations
that have become standard by now. Using the definitions of the $q$-entire 
numbers and the $q$-factorials introduced in the previous section (
\ref{factorials1}, \ref{factorials2} and \ref{factorials3} ), we get 
\begin{equation}
D^N \, \Phi = d^N \, \Phi + {\displaystyle {\sum_{k=1}^{N-1}}} \, 
{N\brack k}_q \, (D^{k-1} \, A) \, d^{N-k}\, \Phi \, + \, (D^{N-1}) \, \Phi
\end{equation}
\indent
The proof uses simple recurrence. Acting again with the operator $D$ and 
separating the two first and the two last terms we get: 
$$D^{N+1} \Phi \, = \, d^{N+1} \Phi + A\,d^N \Phi + {\displaystyle 
{\sum_{k=1}^{N-1}}}\, {N\brack k}_q \, (D^k A)\, d^{N-k} \Phi +  $$
\begin{equation}
+ {\displaystyle {\sum_{k=1}^{N-1}}} \, q^k \, {N\brack k}_q \,(D^{k-1} A) \, 
d^{N-k+1} \Phi \, + \, + q^N \, (D^{N-1} A) \, d\Phi + (D^N A ) \,\Phi
\end{equation}
\indent
Leaving the first and last terms unchanged, and including the term 
$q^N (D^N A)\Phi$ in the second sum, and then shifting the summation index
from $(k-1)$ to $k$, we can re-write the above formula as
$$D^{N+1} \Phi \, = $$
$$\qquad d^{N+1}\Phi + {\displaystyle{\sum_{k=1}^{N+1}}} \,
\Biggl( q^{k+1}\,{N\brack{k+1}}_q \, + \, {N\brack k}_q \, \Biggr) \, (D^k A)
\, d^{N-k} \Phi \,  \, +  (D^N A)\Phi $$
and we get the same formula again if we use the fact that the recurrence
formula (\ref{factorials3}) holds for any value of $q$. In particular, when 
$q$ is a primitive $N$-th root of unity, the formula reduces to
\begin{equation}
D^N \Phi \, = \, (D^{N-1} A) \, \Phi  \, = \, \Omega \, \Phi
\end{equation}
\indent
Let $S$ be an automorphism of the algebra ${\cal{A}}$ in which the $1$-form 
$A$ takes its values. If the module ${\cal{H}}$ is a free one, it induces 
automatically an automorphism of ${\cal{H}}$. It is easy to prove the following
generalization of the ``pure gauge'' connexions and the fact that the
corresponding curvature form must vanish:
\smallskip

\centerline{{\it If  $ A = S^{-1}d S$ and  $ d^N = 0$,   then 
$\Omega = D^{N-1} A = 0$    and  vice-versa.}}
\smallskip

\indent
The proof is by straightforward calculation; indeed, if $A = S^{-1} d S$, then
$ D \Phi = ( d + S^{-1} d S ) \Phi \, = \, d \Phi + S^{-1} d S  \Phi $
(we remind that both $S^{-1}$ and $S$ are of degree $0$ in the sense of 
differential forms);\\
$ D^2 \, \Phi \, = \, (d + S^{-1} dS)^2 \, \Phi = d^2 \Phi + S^{-1} d^2 S \,
\Phi + (1+q) \, S^{-1} dS d\Phi , $\\
$ D^3 \Phi = d^3 \Phi + S^{-1} d^3 S \, \Phi +(1+q+q^2) \, S^{-1} d^2 S \, 
\Phi + (1+q+q^2) S^{-1} dS \, d^2 \Phi; $\\
\indent
and the general formula is
\begin{equation}
D^N \Phi \, = \, d^N \Phi + S^{-1} d^N S \, \Phi + {\displaystyle
{\sum_{i=1}^{N-1}}} \, {n\brack k}_q \, S^{-1} d^k S \, d^{N-k} \Phi
\end{equation}
from which we see that $D^N \Phi = 0$, because for $q$ which is a primitive 
$N$-th root of unity one has $d^N = 0$, and all the symbols ${N\brack k}_q$ 
vanish. This means that $\Omega = D^{N-1} A = 0$, because we have already 
checked that $D^N \Phi = (D^{N-1}A) \, \phi = (\Omega) \, \Phi $
\newline
\indent
The above formulae generalize the notions of connection and curvature in 
a universal way, independent of the realization. In order to make these
formulae useful, one must express the curvature in a more explicit manner,
which shall depend on the realization chosen.
In the next section we show several examples of such realizations, along with
the explicit calculations of the generalized curvature forms.

\section{ Matrix and functional realizations of $q$-differential calculus}

As was stated above, when $q$ is a primitive $n$-th root of unity, the 
$q$-differential algebra can be faithfully represented by $n \times n$ complex
matrices. Most of the results obtained in the previous sections are independent
of realization we choose; nevertheless, some of them can be given explicitly
and lead to the formulae specific for each particular realization. Thus, in
the matrix realization, the ``pure gauge'' connection $1$-form $A$ and the
matrix $\eta$ inducing the exterior ${\Bbb Z}_N$-graded differential can be 
chosen as follows:
\begin{equation}
A = \pmatrix{0& \alpha &0&...&0 \cr 0&0& \beta &..&0 \cr 0&..&..&..&0 
\cr 0&..&..&..& \phi \cr \omega &0&0&..&0}, \, \ \ \,  \eta = 
\pmatrix{0& 1 &0&...&0 \cr 0&0& 1 &..&0 \cr 0&..&..&..&.0
\cr 0&..&..&..& 1 \cr 1 &0&0&..&0},
\end{equation}
whereas the {\it grading matrix} $U$ (which is the $0$-degree generator of
our algebra) is chosen then as
$$U = diag (q, q^2, q^3, ... q^N)$$
The matrix $\Omega = D^{N-1} \, A$ is of degree $N$, i.e. it is diagonal;
moreover, it is easy to check that it is proportional to the unit $N \times N$
matrix with the coefficient $(\alpha +1)(\beta +1)( \gamma +1)...(\omega +1)$. 
Now, taking an arbitrary $0$-degree matrix $S = diag (a, b, c,..., y, z)$, 
its inverse being 
$S^{-1} = diag (a^{-1}, b^{-1}, c^{-1},..., y^{-1}, z^{-1})$, and identifying 
$S^{-1} \, d S $ with $A$, we get
$$ \alpha = (\frac{b}{a} - 1), \ \ \, \beta = (\frac{c}{b} - 1),..., \omega
= (\frac{a}{z} - 1) $$
$$ {\rm {therefore}} \, \ \ \, \ \  (\alpha + 1)(\beta + 1) ... (\omega + 1) 
= (\frac{b}{a})(\frac{c}{b}) (\frac{d}{c})...(\frac{y}{z})(\frac{z}{a}) \, 
= \, 1 $$
and it becomes clear that there is a one-to-one correspondence between the
``pure gauge'' $S^{-1}\,d S$ and null-curvature connections.
\newline
\indent
This realization of the $d^N = 0$ differential calculus seems to be quite
trivial. Hopefully, a more sophisticated version is at hand. Indeed, is it
easy to see that for a given $N$, one can introduce $N$ degree $1$ 
linearly independent complex matrices $\eta_k$, $k = 1,2,..,N$ that can be 
chosen as follows:
\begin{equation}
\eta_k = \pmatrix{0&q^k&0&...&0 \cr 0&0&q^{2k}&...&0 \cr ..&..&..
&...&.. \cr 0 &0&..&0&q^{(N-1)k} \cr
1 &0&0&...&0 }
\end{equation}
Direct calculus shows that
$$ \eta_k^N = (q^k q^{2k} \dots q^{(N-1)k} ) \, {\bf 1}_{N \times N} =
q^{\frac{N(N-1)}{2} k} \, {\bf 1}_{N \times N} $$
which amounts to
$$ \eta_k^N = {\bf 1}_{N \times N} \, \ \ \, {\rm if \ \ N \ \ is \ \ odd 
\ \ or \ \ if \ \ k \ \ is \ \ even} ,$$
\begin{equation}
\eta_k^N = - {\bf 1}_{N \times N} \, \ \ \, {\rm if \ \ N \ \ is \ \ even
\ \ and \ \ k \ \ is \ \ odd}.
\end{equation}
(A similar, though slightly more complicated formula is valid for a totally
symmetrized product of any $N$ generators, which is always proportional to
the unit $N \times N$ matrix, but which in certain cases may vanish).
Now we can define $N$ independent $q$-derivations denoted $d_k, \, k = 1,2,
\dots , N$:
\begin{equation}
d_k \, B = [ \eta_k , B ]_q = \eta_k \, B - q^{{\rm {deg}}(B)} \, B \, \eta_k
\end{equation}
Obviously, the $N$-th power of each of these $q$-differentials vanishes:
$$ d^N_k \, B = 0 \, \ \ \, {\rm for \ \ any \ \ B \in {\cal{A}} } ; $$
but we have also the following generalization of this fact:
\begin{equation}
\displaystyle{\sum_{\pi_{\alpha} (k_1k_2..k_N)}} \, d_{\pi_{\alpha} (k_1)} \,
d_{\pi_{\alpha} (k_2)} \dots d_{\pi_{\alpha} (k_N)} \, B = 0, \, \ \ \, \ \
k = 1, 2, \dots , N.
\end{equation}
where $\displaystyle{\sum_{\pi_{\alpha} (k_1 k_2 .. k_N)}}$ means the sum over 
{\it all} the $N!$ permutations of $N$ indices. In the case when some of the
indices are repeated, one still has:
\begin{equation}
\displaystyle{\sum_{symmetrized}} \, d_{( k_1}^{m_1} \,
d_{k_2}^{m_2} \dots d_{k_p)}^{m_p} \, B = 0 \, \ \ \, \ \ \, \ \ 
{\rm if} \, \ \ \, \ \ m_1 + m_2 + ... + m_p = N
\end{equation}
\indent
At this point we can introduce a generalization of the covariant differential:
it is enough to introduce $N$ independent degree $1$ matrices (``one-forms'')
$A_k, k = 1,2,\dots N$ and to define
$$D_k \, \Phi = d_k \, \Phi + A_k \, \Phi$$
\indent
As in the case of a single $q$-differential $d_q$ , we can show that the
sum over all permutations of $N$ consecutive applications of the covariant
differentials $D_k$ leads to an automorphism of the free module ${\cal{H}}$:
$$\displaystyle{\sum_{all \ \ permutations}} \, D_1 D_2 \dots D_N \, B =
\Omega_{1 2 \dots N} \,  B $$
Furthermore it can be proved by direct calculus that the totally symmetric
quantity $\Omega_{1 2 \dots N}$ is proportional to a unit $N \times N$ matrix,
and it is equal to
\begin{equation}
\Omega_{1 2 \dots (N-1) } \,  = \displaystyle{\sum_{all \ \ permutations}} \, 
( D_1 D_2 \dots D_{N-1} ) \, A_N 
\end{equation}
and similarly in the case of the repeated indices.
\newline
\indent
As an illustration, let us show how this scheme works in the simplest case
when $N = 2$. Then we set
$$ \eta_1 = \pmatrix{0&-1 \cr 1&0} \, , \ \ \, \ \ \eta_2 = \pmatrix{0&1 \cr
1&0} \, , \, \ \ \, A_k = \pmatrix{0&\alpha_k \cr \beta_k &0} \, , \, \ \
k= 1,2. $$
Direct calculus shows that all the components of the curvature $2$-form are
proportional to the $2 \times 2$ unit matrix:
$$\Omega_{12} = \Omega_{21} = (\alpha_1 + \alpha_2 +\beta_1 - \beta_2 +
\alpha_1 \beta_2 + \alpha_2 \beta_1) \, {\bf 1} ; $$
\begin{equation}
\Omega_{11} = ( \alpha_1 - \beta_2 ) \, {\bf 1} , \, \ \ \, \ \ \, \ \ 
\Omega_{22} = ( \alpha_2 + \beta_2 ) \, {\bf 1} 
\end{equation}
If $U$ is a diagonal (degree $0$) matrix, the gauge principle applies now to
the components of the connection $A_k$ as in the classical gauge theory, i.e.
the curvature $N$-form transforms homogeneously when the connection $1$-form
undergoes a gauge transformation:
$$ A_k \rightarrow  \tilde{A_k} = U^{-1} A U + U^{-1} d_k U , \, \ \ \,  
\ \ \, \Omega \rightarrow \tilde{\Omega} = U^{-1} \Omega U . $$
\indent
In the matrix realization the curvature $\Omega$ is always proportional to 
the unit matrix and commutes with all the elements of the algebra, and is 
{\it invariant} undet gauge transformations: $\tilde{\Omega} = \Omega .$
\newline
\indent
An interesting problem arises if we consider the notion of connection forms
with given symmetry properties. The transformation matrices $U$ belong to
the group $GL(N, {\bf R})$. The elements of this group can also act directly
on the set of $N$ matrices $A_k$ inducing a linear transfromation:
$$A_k \rightarrow \tilde{A}_k = M^j_k \, A_j $$
A connection $1$-form is said to be {\it symmetric} with respect to a given
subgroup of the $GL(N, {\bf R})$ group if for any given element of this
subgroup identified as linear transformation $M_k^j$ a corresponding
representation $U(M)$ exists, satisfying
$$ M^j_k \, A_j = U^{-1} A_k U + U^{-1} d_k U $$
Such connections have been considered in the classical case by N. Manton 
(\cite{Manton}) and used in the analysis of Higgs mechanism in electroweak
theory (\cite{Bertrand}).
\newline
\indent
Now, in the usual ${\Bbb Z}_2$-graded case, the most important realization of 
the Grassmann algebra was the {\it algebra of exterior differential forms}
defined on a differential manifold. Here a similar realization can 
be conceived in a ${\Bbb Z}_N$-graded case starting with more explicit calculus
on the example of $N=3$, which has been introduced already in (\cite{Ker2}).
\newline
\indent
We postulate that by definition the differential $d f$ of a function $f$ 
coincides with the usual one:
\begin{equation}
d f = \frac{\partial \,f}{\partial \xi^k} \,d \xi^k =(\partial_k f) \,d\xi^k
\end{equation}
\indent
When computing formally higher-order differentials, we shall suppose that
our exterior differential operator $d$ obeys the generalized {\it graded}
Leibniz rule:
\begin{equation}
d\,(\omega \, \phi) = d\,\omega\,\phi + q^{{\rm{deg}}(\omega)}\,\omega\,d\,\phi
\end{equation}
where we suppose that $q$ is an $N$-th order root of unity, instead of $-1$ 
in the ${\Bbb Z}_2$-graded case, and that the grades add up modulo $N$ under 
the associative multiplication of exterior forms; the functions are of grade 
$0$, and the operator $d$ raises the grade of any form by $1$, which means 
that the linear operator $d$ applied to $\xi^k$ produces a 1-form whose 
${\Bbb Z}_N$-grade is $1$ by definition; when applied two times, by iteration,
it will produce a new entity, which we shall call a {\it 1-form of grade $2$,} 
denoted by $d^2 \xi^k$. Finally, we require that $d^N = 0$. 
\newline
\indent
Let {\it F} denote the algebra of functions $C^{\infty} (\xi^k) $, over 
which the ${\Bbb Z}_N$-graded algebra generated by the forms $d \xi^i$, 
$d^2 \xi^k$, $d^3 \xi^k$, etc., behaves as a {\it left} module. 
In other words, we shall be able to multiply the forms $d\xi^i$ , $d^2 \xi^k$,
 $d\xi^i d\xi^k$ , etc. 
by the functions {\it on the left} only; right multiplication will just not 
be considered here. That is why we shall write by definition, e.g.                 
\begin{equation}
d (\xi^i \xi^k) := \xi^i d\xi^k + \xi^k d\xi^i
\end{equation}
\indent
This amounts to suppose that the coordinates (functions) commute with the
$1$-forms, but do not necessarily commute with the forms of higher order.
With thus established ${\Bbb Z}_3$-graded Leibniz rule, the postulate $d^3 = 0$ 
suggests in an almost unique way the ternary and binary commutation rules
for the differentials $d\xi^i$ and $d^2 \xi^k$. Consider the
differentials of a function of the coordinates $\xi^k$:
\begin{eqnarray*}
df :&=&  (\partial_i f) d\xi^i \ \ ; \ \ \  \ d^2 f := (\partial_k 
\partial_i f) d\xi^k d\xi^i + (\partial_i f) d^2 \xi^i \ ; \nonumber\\ 
d^3 f &=& (\partial_m \partial_k \partial_i f) d\xi^m d\xi^k d\xi^i + 
(\partial_k \partial_i f) d^2 \xi^k d\xi^i \nonumber\\
&&{}+ j (\partial_k \partial_i f) 
d\xi^i d^2 \xi^k + (\partial_k \partial_i f) d\xi^k d^2 \xi^i \ ;\nonumber 
\end{eqnarray*}
(we remind that the last part of the differential, $(\partial_i f)d^3 \xi^i$, 
vanishes by virtue of the postulate $d^3 \xi^i = 0$). Supposing that partial 
derivatives commute, exchanging the summation indices $i$ et $k$ in the last 
expression and replacing $1+j$ by $-j^2$, we arrive at the following two 
conditions that lead to the vanishing of $d^3 f$ :
\begin{equation}
d\xi^m d\xi^k d\xi^i + d\xi^k d\xi^i d\xi^m + d\xi^i d\xi^m d\xi^k = 0 \ ; 
\  d^2 \xi^k d \xi^i - j^2  d \xi^i d^2 \xi^k = 0.
\end{equation}
\indent
which lead in turn to the following choice of relations:
\begin{equation}
d \xi^i d \xi^k d \xi^m = j d \xi^k d \xi^m d \xi^i , \  \ {\rm and} \  \  
d \xi^i d^2 \xi^k = j d^2 \xi^k d \xi^i .
\end{equation}
\indent
Strictly speaking, the above formulae hold only for the {\it symmetric} part
of the above expression; we choose to impose stronger relations in order to
make the resulting space of forms finite-dimensional.
\newline
\indent
Extending these rules to {\it all} the expressions with a well-defined grade, 
and applying the associativity of the ${\Bbb Z}_3$-exterior product, we see 
that all products of the type $d \xi^i d \xi^k d \xi^m d \xi^n$ and 
$d \xi^i d\xi^k d^2\xi^m$ must vanish, and along with them, also the 
monomials of higher order containing these as factors.
\newline
\indent
Still, this is not sufficient in order to satisfy the rule $d^3 = 0$ on all 
the forms spanned by the generators $d \xi^1$ and $d^2 \xi^k$. It can be 
proved without much pain that the expressions containing $d^2 \xi^i d^2 \xi^k$ 
must vanish, too; so we set forward the additional rule declaring that any 
expression containing {\it four or more} operators $d$ must identically 
vanish. With this set of rules we can check that $d^3 = 0$ on all the forms, 
whatever their grade or degree. 
\newline
\indent
As in the case of the matrix algebra realization, it is very easy to introduce 
a ${\Bbb Z}_N$-graded generalization. Replacing $j$ by $N$-th primitive root of 
unity, $q$, and introducing higher "multi-differentials", $d^3 \, \xi^k$,
$d^4 \, \xi^k$,..., up to $d^{N-1} \, \xi^k$ , we compute next differentials
as follows: we impose on the operator $d$ the $q$-graded Leibniz rule and
we require that $d^N = 0$, we can impose the following minimal set of 
generalized commutation rules on the products of forms of order $N$:
\begin{equation}
d \xi^{k_1} d \xi^{k_2} \dots d\xi^{k_N} = q\, d \xi^{k_2} \dots d \xi^{k_N} 
d \xi^{k_1} = q^2 \, d \xi^{k_3}  \dots
d \xi^{k_N} d\xi^{k_1} d \xi^{k_2} ... {\rm etc.,}
\end{equation}
As a corollary, one can conjecture that for $N \geq 3$ any product of more 
than $N$ such $1$-forms must vanish, which has been proved for the general 
${\Bbb Z}_N$-graded Grassmann algebras.

\indent
As now $d^2 \neq 0$, $d^3 \neq 0$,..., $d^{N-1} \neq 0$, we must introduce
new independent differentials, $ d^2 \xi^k , \ \  d^3 \xi^k , \ \ \dots , \ \ 
d^{N-1} \xi^k .$ Each kind of these new {\it ``one-forms of degree $m$''}, 
with $m = 1,2 \dots (N-1)$ spans a basis of a $D$-dimensional linear space.
\newline
\indent
We assume that all the products of forms whose total degree is less than $N$ 
are independent and span new modules over the algebra of functions with 
appropriate dimensions, e.g. the products of degree $2$, $d\xi^k d \xi^m$, 
span a $D^2$-dimensional linear space; so do the products $d^2 \xi^k d \xi^m$ 
and, independently, $d \xi^m d^2 \xi^k$ (if $D > 3$), and so on. On the other 
hand, all other products of degree $N$ must obey the following commutation 
relations, which are compatible with the cyclic commutation relations for 
the product of $N$ $1$-forms, for example:
$$d^p \xi^k d^{N-p} \xi^l = q^p \, d^{N-p} \xi^l d^p \xi^k, $$
\begin{equation}
d^{N-p} \xi^k d \xi^{l_1} d \xi^{l_2} \dots d \xi^{l_p} = q^{N-p} \, 
d \xi^{l_1} d \xi^{l_2} \dots d \xi^{l-p} d^{N-p} \xi^k, \, {\rm etc.}
\end{equation}
\indent
Finally, we shall assume that not only the products of $N+1$ and more $1$-forms
vanish, but along with them, also any other products of all kinds of forms
whose total degree is greater than $N$. This additional assumption is necessary
in order to ensure the coordinate-independent character of the condition
$d^N = 0$, because under a coordinate change all the products of forms of 
given order mix up and transform into each other, e.g. the terms
$d \xi^j d \xi^k d \xi^l$ with the terms of the type $d^2 \xi^k d \xi^l$, and 
similarly for higher order terms. 
\newline
\indent
It is easy to prove that for a given $N$ it is enough to assume $d^N \xi^k = 0$
and the $N$-cyclic commutation rule 
$$d \xi^{k_1} d \xi^{k_2} d \xi^{k_3} \dots d \xi^{k_N} = q\, d \xi^{k_2}
d \xi^{k_3} \dots d \xi^{k_N} d \xi^{k_1}$$
implemented with its generalization for any product of two exterior forms
of the total order adding up to $N$,
$$\omega \, \phi = q^{p(N-p)} \phi \, \omega = q^{-p^2} \phi \, \omega $$
whenever deg$(\omega) = p$ and deg$(\phi) = N - p$, in order to ensure 
that $d^N f = 0$, and in general, $d^N \omega = 0$ for any differential form 
$\omega$.
\newline
\indent
We end this section by showing how the gauge-invariant analogs of curvature
are expressed via the unique gauge-invariant curvature $2$-form $F_{ik}$ and 
its covariant derivatives. Indeed, for $N = 3$ we have (cf. Ref \cite{Ker2}):
\begin{equation}
\Omega = D^2 A = \frac{1}{3}\, (D_i F_{km} + j D_m F_{ki}) \, d \xi^i d \xi^k
d \xi^m + F_{ik} \, d^2 \xi^i d \xi^k  \, ,
\end{equation}
whereas when $N = 4$ we get, with 
$\Theta_{iklm} = \frac{1}{4} (D_i D_k F_{lm} + iD_m D_l F_{ki})$,\ 
$ \Phi_{klm} = D_k F_{lm} + i D_m F_{kl}$, \ 
$ \tilde{F}_{ik} = \frac{i}{2} \, F_{ik}$:
\begin{eqnarray*}
\Omega &=& D^3 A \nonumber\\
&=& \Theta_{jklm} d \xi^j d \xi^k d \xi^l d \xi^m +
\Phi_{ikm} d^2 \xi^i d \xi^k d \xi^m + F_{ik} d^3 \xi^i d \xi^k +
\tilde{F}_{ik} d^2 \xi^i d^2 \xi^k.
\end{eqnarray*}

\section{ Deformations of the covariant $q$-differential}

The $q$-differentials defined on the ${\Bbb Z}$-graded associative unital 
algebra can be embedded in a bundle over the complex plane ${\Bbb C}^1$, with 
a typical fibre being the space of linear operators defined on the algebra 
{\cal{A}}.
\newline
\indent
With this view on the $q$-differentials, a natural question can be asked now, 
namely, what would be the effect of one or more {\it deformations} of the 
complex parameter $q$ itself, $q \rightarrow q + \epsilon$, $\epsilon 
\, \in {\Bbb C}^1$ ? In particular, it is interesting to see what will happen 
if we perform a series of such deformations, returning back to the initial
value of $q$. In other words, if we consider a series of infinitesimal 
complex deformations $\epsilon, \, {\epsilon}^{'}, \, {\epsilon}^{''}, \dots , 
{\epsilon}^{(n)}$, such that $\epsilon + {\epsilon}^{'} + {\epsilon}^{''} + 
\dots + \epsilon^{(n)} \, = 0$, what will be the result of consecutive action 
of $n$ corresponding $(q+ \epsilon^{(k)})$-deformed differentials, or of their
combinations with various permutations ? This question seems quite pertinent
in the context of $q$-differential calculus, and nearly as natural as the
investigation of consecutive infinitesimal diffeomorphisms of a manifold,
induced by two different vector fields, that led Sophus Lie to the definition
of the Lie bracket as the commutator of two such automorphisms.
\newline
\indent
Let us consider first the simplest case of $N = 2$. According to the scheme
exposed above, we shall put ${\epsilon}^{'} = - \epsilon$. Now, let us compare
the consecutive action of the operators $d_{q+\epsilon}$ and $d_{q-\epsilon}$
taken in two different orders. Acting on an arbitrary element $B$ of degree
$1$, for example, we have by definition:
$$d_{q+\epsilon} \, B = \eta \, B - (q+\epsilon) \, B \, \eta ;$$
and taking into account that now $d_{q+\epsilon}\,B$ is of degree $2$, i.e.
$0_{mod \, 2}$, we arrive at the following expression:
$$d_{q-\epsilon} d_{q+\epsilon} \, B = \epsilon (B \eta^2 - \eta B \eta) =
\epsilon (B - \eta B \eta) \, \ \ \, (q= -1 ,  \eta^2 = {\bf 1}) $$
\indent
Evaluating similar expression with $d_{q+\epsilon} \, d_{q-\epsilon}$ taken 
in the reverse order amounts to changing the sign of the parameter $\epsilon$ 
in the above expression. Therefore adding or substracting the two expressions 
we get
\begin{eqnarray}
\biggl[ \,d_{q-\epsilon} \, d_{q+\epsilon} \, + \, d_{q+\epsilon} \,
d_{q-\epsilon} \, \biggr] \, B &=& 0 ,\\
\biggl[ \,d_{q-\epsilon} \, d_{q+\epsilon} \, - \, d_{q+\epsilon} \,
d_{q-\epsilon} \, \biggr] \, B \, &=& 2 \epsilon (B - \eta B \eta)
\end{eqnarray}
\indent
When the degree of B is $0$, the corresponding formula reads: 
\begin{eqnarray}
\biggl[ \,d_{q-\epsilon} \, d_{q+\epsilon} \, + \, d_{q+\epsilon} \,
d_{q-\epsilon} \, \biggr] \, B &=& 0 , \\
\biggl[ \,d_{q-\epsilon} \, d_{q+\epsilon} \, - \, d_{q+\epsilon} \,
d_{q-\epsilon} \, \biggr] \, B \, &=& 2 \epsilon (B + \eta B \eta)
\end{eqnarray}
It is interesting to note that the operators defined as \\
\centerline{${\cal{P}}_1 (B) =\frac{1}{2}(B - \eta B \eta)$ and 
${\cal{P}}_2 (B) = \frac{1}{2}(B+ \eta B \eta)$}\\ are projectors, 
i.e. ${\cal{P}}_k^2 (B) = 
{\cal{P}}_k (B)$, $k=1,2.$ Obviously, ${\cal{P}}_1 \, {\cal{P}}_2 =
{\cal{P}}_2 \, {\cal{P}}_1 = 0$, so that these operators project onto the
subspaces of degree $1$ and $0$ respectively.
\newline
\indent
Similar result can be obtained for a series of 3 consecutive deformed
$q$-differentials around the value of $q = j = e^{2 \pi i /3}$. Let us note
$$d_1 \, = \, d_{q+j \epsilon}, \, \ \ \, d_2 \, = \, d_{q+ j^2 \epsilon}
\, \ \ {\rm \, and \,} \, \ \ \, d_3 = d_{q+\epsilon} . $$
\indent
Direct calculus leads to the following :
\begin{equation}
[ \, d_1 d_2 d_3 + d_2 d_3 d_1 + d_3 d_1 d_2 \, ] \, B = 0
\end{equation}
whereas the combination that generalizes the {\it commutator} for the 
${\Bbb Z}_3$-graded case produces (up to a factor) projection operators
$$ \frac{1}{9} [ \, d_1 d_2 d_3 + j d_2 d_3 d_1 + j^2 d_3 d_1 d_2 \, ] \, B =
\frac{1}{3} (B + j \eta B \eta^2 + j^2 \eta^2 B \eta) = P_1 (B).$$
if deg$(B) = 1$, and a similar operator $P_2$ with $j$ and $j^2$ interchanged 
when deg$(B) = 2$. Here again, $P_1^2 = P_1, \, \ \ P_2^2 = P_2, \, \ \ 
P_1 P_2 = P_2 P_1 = 0$.
\newline
\indent
Similar projection operators may be obtained by applying the same scheme to
the properly symmetrized products of $N$ deformed differentials in the case 
of a ${\Bbb Z}_N$-graded generalization. More complicated scheme arises when we
consider all the $N$ independent ``partial'' $q$-differentials.
\newline
\indent
A natural and important question to ask is the following : how the covariant
$q$-differentials react under these cyclic deformations ? What kind of an
operator we get applying the $N$-fold products of covariant $q$-differentials 
deformed as follows:
\begin{equation}
D_q = d_q + A \Rightarrow \tilde{D}_{q+\epsilon} = d_{q+\epsilon} + A + 
\epsilon \Lambda
\end{equation}
In the $N=2$ matrix realization case, the symmetric part of the product yields
the usual curvature and no extra terms linear in the deformation parameter
$\epsilon$, while the anti-symmetrized product is proportional to $\epsilon$
and is equal to :
$$ [ \, D_{q+\epsilon} \, {\tilde{D}}_{q-\epsilon} - D_{q-\epsilon} \, 
{\tilde{D}}_{q+\epsilon} \, ] \, B = D_q \, B \eta + (D_q \Lambda) B $$
\indent
In conclusion, we think that the last idea consisting in a discrete analogue
of the holonomy, might be useful for the analysis of new symmetries of the
$q$-deformed covariant differential operators and for the classification of
their invariant properties. In some sense this is an analog of the geodesic
deviation equation in classical differential geometry.

\acknowledgements{ The authors gratefully acknowledge many helpful 
remarks and hints by M. Dubois-Violette and enlightening discussions with 
B. Le Roy.}

\end{document}